\documentclass[fleqn]{mat01}
\usepackage{times,mathtimy,amssymb}
\begin{document}

\setcounter{page}{15}
\firstpage{15}

\font\zz=msam10 at 10pt
\def\Box{\mbox{\zz{\char'244}}}

\font\www=mtgub at 10.4pt
\def\piy{\mbox{\www{\char'160}}}
\def\Delt{\mbox{\www{\char'104}}}

\font\bi=tibi at 10.4pt

\newtheorem{theore}{Theorem}
\renewcommand\thetheore{\arabic{section}.\arabic{theore}}
\newtheorem{theor}[theore]{\bf Theorem}
\newtheorem{propo}[theore]{\rm PROPOSITION}
\newtheorem{lem}[theore]{Lemma}
\newtheorem{definit}[theore]{\rm DEFINITION}
\newtheorem{coro}[theore]{\rm COROLLARY}
\newtheorem{rem}[theore]{Remark}
\newtheorem{exampl}[theore]{Example}

\newtheorem{proo}{Proof}
\renewcommand\theproo{(\arabic{proo})}
\newtheorem{prof}[proo]{Proof of}
\newtheorem{sttp}{Step}
\renewcommand\thesttp{\arabic{sttp}}
\newtheorem{step}[sttp]{Step}
\def\claim{\trivlist\item[\hskip\labelsep{{\it Claim.}}]}
\def\potc{\trivlist\item[\hskip\labelsep{{\it Proof of the Claim.}}]}

\newtheorem{ppro}{Proof}
\renewcommand\theppro{(\arabic{section}) $\Rightarrow$ (\arabic{ppro})}
\newtheorem{proff}[ppro]{Proof of}

\newcommand{\cc}{{\mathcal C}}
\newcommand{\cq}{{\mathcal Q}}
\newcommand{\cs}{{\mathcal S}}
\newcommand{\ca}{{\mathcal A}}
\newcommand{\cp}{{\mathcal P}}

\newcommand{\br}{{\mathbb R}}
\newcommand{\bc}{{\mathbb C}}
\newcommand{\bz}{{\mathbb Z}}
\newcommand{\bs}{{\mathbb S}}
\newcommand{\bp}{{\mathbb P}}

\newcommand{\lr}{\longrightarrow}
\newcommand{\ra}{\rightarrow}
\newcommand{\hra}{\hookrightarrow}

\newcommand{\ub}{\underline{b}}
\newcommand{\ut}{\underline{t}}
\newcommand{\uc}{\underline{c}}

\renewcommand{\phi}{\varphi}
\renewcommand{\bar}{\overline}

\renewcommand{\theequation}{\thesection\arabic{equation}}

\title{On the fundamental group of real toric varieties}

\markboth{V~Uma}{On the fundamental group of real toric varieties}

\author{V~UMA}

\address{Institute of Mathematical Sciences, C.I.T. Campus,
Chennai~600~113, India\\
\noindent E-mail: uma@imsc.res.in}

\volume{114}

\mon{February}

\parts{1}

\Date{MS received 25 March 2003; revised 10 October 2003}

\begin{abstract}
Let $X(\Delta)$ be the real toric variety associated
to a smooth fan $\Delta$. The main purpose of this article is: (i) to
determine the fundamental group and the universal cover of
$X(\Delta)$, (ii) to give necessary and sufficient conditions on
$\Delta$ under which $\pi_1(X(\Delta))$ is abelian, (iii) to give
necessary and sufficient conditions on $\Delta$ under which
$X(\Delta)$ is aspherical, and when $\Delta$ is complete, (iv) to give
necessary and sufficient conditions for $\cc_{\Delta}$ to be a
$K(\pi,1)$ space where $\cc_{\Delta}$ is the complement of a real
subspace arrangement associated to $\Delta$.
\end{abstract}

\keyword{\looseness 1 Real toric varieties; fundamental group; asphericity;
$K(\pi,1)$ subspace\break arrangements.}

\maketitle

\setcounter{section}{-1}
\section{Notations}\vspace{.2pc}

\noindent $N\cong\bz^n, M = {\rm Hom} (N,\bz)$ and $\langle,\rangle= $
the dual pairing.\\
$N_\br = N \bigotimes_{\bz}\br$. $\Delta = $ smooth fan in $N_\br
;\sigma$ and $\tau$ denote cones in $\Delta$.\\
Let $\sigma$ be a cone in $\Delta$. $S_{\sigma}=\sigma^{\vee}\cap
M = \{ u \in M:\langle u,v \rangle \geq 0 ~ \forall~ v \in \sigma \}$.\\
$\Delta(k)$ = cones of dimension $k$. $\Delta(1)$ are the edges
and $\#\Delta(1)=d$.\\
$\Delta(1)= \{\rho_1, \rho_2, \ldots, \rho_d\}$. Let $v_j$ be the
primitive vector along $\rho_{j}$ then, $\langle
v_{i_1},\ldots,v_{i_k}\rangle$ denotes the cone spanned by
$\{v_{i_1},\ldots,v_{i_k}\}$.\\
$(U_\sigma)_{\bc} = {\rm Hom}_{sg}(S_\sigma,\bc), U_\sigma ={\rm
Hom}_{sg}(S_\sigma,\br)$ and $(U_\sigma)_{+}={\rm
Hom}_{sg}(S_\sigma,\br_{+}) \forall~ \sigma\in\Delta$ where $\br_{+} =
\br^{+} \cup \{ 0 \}$. Here, ${\rm Hom}_{sg}$ denotes the semigroup
homomorphisms which sends $0$ to $1$.\\
$X = $ smooth real toric variety of dimension $n$ associated to
$\Delta$.\\
$X_{\bc} = $ the complex toric variety whose real part is $X$.\\
$X_{+} = $ the non-negative part of $X$.\\
$T_2:= {\rm Hom} (M,\bz_2)\hookrightarrow T_{\br} := U_{\{0\}}=
{\rm Hom} (M,\br^*); T_{\bc}:= {\rm Hom} (M,\bc^*); T_{+}:= {\rm Hom}
(M,\br^{+})$.\\
For every $\sigma \in \Delta, x_\sigma\in U_\sigma$ is a {\it
distinguished point} defined as:
\begin{equation*}
x_\sigma(u) = \begin{cases}
1 & \forall~ u \in \sigma^{\perp}\\
0 &{\rm otherwise} \end{cases}.
\end{equation*}
$O_\tau = $ orbit of $x_\tau$ under the action of
$T^{'}\simeq(\br^*)^n$ and $V(\tau)= \overline{O_\tau}$.\\
${\rm Stab}(x_{\tau})= $ stabilizer of $x_{\tau}$ under the action
of $T_{2}$.\\
$(O_\tau)_{+} = $ orbit of $x_\tau$ under the action of
$(\br^{+})^n$ and $V(\tau)_{+}= \overline{(O_\tau)}_{+}$.\\
$W(\Delta) =\langle~ s_j:j=1,2,\ldots,d~ \mid~ s_j^{2}: 1\leq j \leq
d,~ (s_is_j)^2$ whenever $\langle v_i,v_j\rangle\in \Delta~\rangle$.
Then, $W(\Delta)$ is a {\it right-angled Coxeter group} associated to
$\Delta$. In many places when the context is clear, we shall denote
$W(\Delta)$ simply by $W$.\\
$S_{N}:=(N_{\br}-\{0\})/\br_{>0}$ be the sphere in $N_{\br}$ and
let $\pi:N_\br-\{0\}\lr S_{N}$ be the projection.

\section{Introduction}

Let $\Delta$ be a smooth fan in the lattice $N\cong\bz^n$. Let
$X(\Delta)_{\bc}$ be the complex toric variety associated to $\Delta$.
Let $X(\Delta)_\br$ be the real part of $X(\Delta)_{\bc}$ which we
call the real toric variety associated to $\Delta$. We shall denote
$X(\Delta)_\br$ by $X(\Delta)$ for convenience, as it is going to be
our main object of study. For the definition and basic facts on real
toric varieties (cf. ch.~4 of \cite{Fult} and \S2 of \cite{Jurk}).
We mainly follow \cite{Fult} for notations and background material on
toric varieties.

In this paper we describe the fundamental group and the universal cover
of $X(\Delta)$. We were motivated by the paper \cite{DJ} of Davis and
Januszkiewicz (cf. Cor.~4.5, p.~415 of \cite{DJ}), where they prove the
corresponding results for the real part of a toric manifold (now also
known as a quasitoric manifold). We show that the same results can be
obtained for a real toric variety $X(\Delta)$ associated to a smooth fan
$\Delta$ not necessarily complete, the basic tool being the theory of
{\em developments of complexes of groups} in chapter~II.12 of
\cite{Brid}. We further give necessary and sufficient conditions on
$\Delta$ for $X(\Delta)$ to be aspherical, motivated by the recent
papers of Davis, Januszkiewicz and Scott (cf. Theorem~2.2.5, p.~27 of
\cite{DJS}), where they again prove similar results for a small-cover.
For this purpose, we too rely primarily on the results of Davis
\cite{Davis}, however in many places we give different proofs using the
technique of {\em development} which is consistent with the theme of
this paper (cf. Theorem~\ref{aspher1} in \S{6}).

Besides generalizing the previous results to the setting of a smooth
real toric variety $X(\Delta)$, we give a presentation for the
fundamental group $\pi_1(X(\Delta))$ completely in terms of the fan;
furthermore, we give necessary and sufficient conditions on $\Delta$
under which $\pi_1(X(\Delta))$ is abelian, and we also show that the
torsion elements are always of order $2$.

Finally in \S7, when $\Delta$ is complete, we relate $X(\Delta)$ to
$\cc_{\Delta}$ which is the complement of a real coordinate subspace
arrangement in $\br^d$ where $d$ is the number of edges in $\Delta$. We
could call $\cc_{\Delta}$ {\em the real toric subspace arrangement
associated to $\Delta$}. It is nothing but the real analogue of
the complement of the complex subspace arrangement in $\bc^d$, for which
$X(\Delta)_{\bc}$ is realized as the geometric quotient under the action
of $(\bc^*)^{d-n}$ (cf. \cite{Audin,Cox}). Moreover, since
$\cc_{\Delta}$ is homotopically equivalent to a covering space of
$X(\Delta)$, we can describe its fundamental group and give conditions
on $\Delta$ under which it is a $K(\pi,1)$ space.

Finding $K(\pi,1)$ arrangements seems to be an interesting problem in
topology (cf. \cite{Orlik} and \cite{khov}) and we get many such
examples. Similar to the results of \cite{khov}, in our case too it
turns out that $\cc_{\Delta}$ is $K(\pi,1)$ if and only if it is the
complement of a union of precisely codimension 2 subspaces (cf.
\cite{Buck,DJS} for other results related to real subspace
arrangements).

Before we state the main theorems let us fix the following terminology:

\begin{enumerate}
\renewcommand\labelenumi{}
\item Let $\Delta(1)$ denote the edges of $\Delta, d=\#\Delta(1)$, and let
$\{v_1,v_2, \ldots,v_d\}$ denote the primitive vectors along the
edges. We assume that $\{v_1,v_2,\ldots,v_n\}$ form a basis for the
lattice $N$ and let $\{u_1,\ldots,u_n\}$ be the dual basis in $M$.

\item Let $W(\Delta)=\langle s_{j_1},\ldots, s_{j_d} \mid s_j^{2}: 1\leq j
\leq d, (s_i s_j)^2$ whenever $\{v_i,v_j\}$ spans a cone in $\Delta
\rangle$ be the right-angled Coxeter group associated to $\Delta$.

\item We call the fan $\Delta$ {\em flag-like} if and only if the following
condition holds for every collection of primitive edge vectors
$\{v_{i_{1}},\ldots, v_{i_{r}}\}$ in $\Delta$: if for every $1\leq k,
l\leq r $, $\{ v_{i_{k}}, v_{i_{l}}\}$ spans a cone in $\Delta$, then
$\{ v_{i_{1}},\ldots, v_{i_{r}}\}$ together spans a cone in $\Delta$.

\item Let $X(\Delta)$ be a smooth and connected real toric variety.
\end{enumerate}
We now state the main results in the paper.

\begin{theor}[\!]\label{abe} The fundamental group $\pi_{1}(X(\Delta))$
is abelian if and only if one of the following holds in $\Delta$.
\begin{enumerate}
\renewcommand\labelenumi{{\rm (\roman{enumi})}}
\leftskip .2pc
\item For every $1\leq i,j\leq d, \{v_i,v_j\}$ spans a cone in
$\Delta$. In this case{\rm ,} $\pi_1(X)$ is isomorphic to $\bz_2^{d-n}$.

\item For each $~1\leq j\leq d$ there exists at most one $i=i_j$ with
$1\leq i_j\leq n$ such that{\rm ,} $\{v_{i_j},v_j\}$ does not span a cone in
$\Delta$ and $\langle u_{i_j},v_j\rangle =1~{\rm mod}~2$. Further{\rm ,} for each
$n+1\leq k\leq d~$ such that $k \neq j$ we have{\rm ,} $\langle
u_{i_j},v_k\rangle=0~{\rm mod}~2$.\vspace{-.5pc}
\end{enumerate}
\end{theor}

\begin{theor}[\!]\label{asp} The real toric variety $X(\Delta)$ is
aspherical if and only if $\Delta$ is flag-like.
\end{theor}

\begin{theor}[\!]\label{arr} Let $\cc_{\Delta}$ be the complement of the
subspace arrangement related to $\Delta$ as above. Then{\rm ,}
$\pi_1({\mathcal C}_{\Delta})$ is isomorphic to the commutator subgroup
of $W(\Delta)$. Further{\rm ,} $\cc_{\Delta}$ is aspherical if and only
if it is the complement of a union of precisely codimension $2$
subspaces.
\end{theor}

We prove Theorem~\ref{abe} in \S5, Theorem~\ref{asp} in \S6 and
Theorem~\ref{arr} in \S7.

In this context, we also mention that the cohomology ring with $\bz_2$
coefficients of smooth, complete real toric varieties and
$H_1(X(\Delta),\bz)$ has been studied by Jurkiewicz (cf. \cite{Jurk}).

\section{The universal cover of {\bi X}\hbox{{\bf ($\Delta$)}}}

In this section we shall determine the universal cover and the
fundamental group of $X$. For this purpose, we primarily apply the
contents of pp.~367--386 of ch.~II.12 of \cite{Brid}.

We begin with the elementary topological description of a real toric
variety in the following proposition. The proof essentially follows
from the proposition on p.~79, ch.~4 of \cite{Fult} by replacing
$X_{\bc}$ by $X$ and $\bs^{1}$ by $\bs^{1}\cap \br\simeq\bz_2$. For
details, also see p.~36, \S3\break of \cite{Jurk}.

\setcounter{theore}{0}
\begin{propo}\label{RT}{\rm \cite{Fult,Jurk}}$\left.\right.$\vspace{.5pc}

\noindent There is a retraction{\rm ,}
$X_+\stackrel{i}{\hookrightarrow} X\stackrel{r}{\rightarrow} X_{+}$
given by the absolute value map{\rm ,} $x\mapsto |x|$ from $\br_+
\subset \br \rightarrow \br_+$ which identifies $X_+$ with the quotient
space of $X$ by the action of the compact real $2$-torus{\rm ,}
$T_2={\rm Hom}(M,\bz_2)$. Further{\rm ,} there is a canonical mapping $T_2
\times X_+\rightarrow X$ which realizes $X$ as a quotient space{\rm ,}
$T_2\times X_+/\sim$ where{\rm ,} $(t,x)\sim (t',x')$ if and only if
$x=x'$ and $t\cdot(t')^{-1}\in {\rm Stab}(x_{\tau})$ where{\rm ,} $x\in
(O_{\tau})_{+}$. The retraction $X \rightarrow X_+$ maps $O_\tau$ to
$(O_\tau)_+$ and $V(\tau)$ to $V(\tau)_+$ and the fiber over
$(O_\tau)_+$ is $T_{\tau}:={\rm Hom}(\tau^{\perp}\cap M,\bz_2)$ which is a
compact real $2$-torus of dimension $n - {\rm dim}(\tau)$.
\end{propo}

We now observe the following property of $X_{+}$.

\begin{lem}\label{SC} $X_{+}$ is contractible.
\end{lem}

\begin{proof} Recall that $x_{\{0\}}$ is the distinguished point of
$(U_{\{0\}})_{+}\simeq (\br^+)^n$. We first show that for every
$\sigma\in\Delta$, $(U_{\sigma})_{+} = {\rm
Hom}_{sg}(S_{\sigma},\br_{+})$ is contractible to the point
$x_{\{0\}}\in (U_{\{0\}})_{+}\subseteq (U_{\sigma})_{+}$. This is
because, $(1-t)\cdot x+ t\cdot x_{\{0\}}\in {\rm
Hom}_{sg}(S_{\sigma},\br_{+})$ for every $t\in I=[0,1]$. The only thing
we need to check here is that, if both $u,-u\in S_{\sigma}$ then
$((1-t)\cdot x+ t\cdot x_{\{0\}})(u)=(1-t)\cdot x(u)+ t\cdot
x_{\{0\}}(u)>0$. This clearly holds since, $x(u)>0$ and $\br^{+}$ is
convex. Thus the map $H_{\sigma}: (U_{\sigma})_{+}\times I\rightarrow
(U_{\sigma})_+$ defined as $H_{\sigma}(x,t)=(1-t)\cdot x+t\cdot
x_{\{0\}}$ is a strong deformation retraction of $(U_{\sigma})_{+}$ to
the point $x_{\{0\}}$.

Furthermore, by definition, the $H_{\sigma}$'s for $\sigma\in \Delta$ are
compatible with the inclusions $(U_{\tau})_{+}\subseteq
(U_{\sigma})_{+}$ whenever $\tau<\sigma$ in $\Delta$. Therefore, since
$X_{+}$ is the union of $(U_{\sigma})_{+}$'s for $\sigma\in \Delta$, we
can glue together the maps $\{H_{\sigma}\}_{\sigma\in\Delta}$ to get a
strong deformation retraction $H$ of $X_+$ to $x_{\{0\}}$. Hence the
lemma.\hfill $\Box$
\end{proof}

\begin{propo}$\left.\right.$\vspace{.5pc}\label{Coxeter}

\noindent Let $\Delta$ be a smooth fan. We then have the following{\rm :}
\begin{enumerate}
\renewcommand\labelenumi{{\rm \arabic{enumi}.}}
\item $(X_+,(V(\tau)_+)_{\tau\in \Delta})$ is a stratified
space with strata $\{V(\tau)_+\}_{\tau \in \Delta}$ indexed by the
poset ${\Delta}$.

\item Associated to this stratified space we have a simple complex of
groups $G({\Delta})=(G_{\tau},\psi_{\sigma\tau})$ where the local group
at the stratum $V({\tau})_{+}$ is $G_{\tau}={\rm Stab}(x_{\tau})$ under
the action of $T_2={\rm Hom}(M,\bz_{2})$ and $\psi_{\sigma\tau}:G_{\tau}
\rightarrow G_{\sigma}$ {\rm (}for $\tau<\sigma$ in $\Delta${\rm )} are
canonical inclusions and we have a simple morphism
$\varphi=(\varphi_{\tau}): G(\Delta)\rightarrow T_{2}\simeq\bz_{2}^n$
injective at the local groups.

\item For the above simple complex of groups $G{
(\Delta)}=(G_{\tau},\psi_{\sigma\tau})${\rm ,} the direct limit
$\widehat{G(\Delta)}$ is isomorphic to $W(\Delta)$. We therefore have a canonical
simple morphism $\iota=(\iota_{\tau}):G(\Delta)\rightarrow W(\Delta)$.
\end{enumerate}
\end{propo}

\begin{proof}$\left.\right.$

\begin{prof}{\rm
Since the orbit space decomposition of $X_{+}$
under the action of $T_{+}$ is obtained by restriction of scalars from
that of $X_{\bc}$ under the action of $T_{\bc}$, it follows that
$(X_+,V(\tau)_+)$ is a stratified space with strata $V(\tau)_+$
indexed by $\Delta$.}
\end{prof}

\begin{prof}{\rm Let $G_{\tau}={\rm Stab}(x_{\tau})\subseteq T_2$. We
then have canonical inclusions, $\psi_{\sigma\tau}:G_{\tau}\hra
G_{\sigma}$ whenever $\tau<\sigma$ in $\Delta$ and,
$\varphi_{\tau}:G_{\tau}\hra T_2$ for every $\tau\in\Delta$. Then
$G(\Delta)=(G_{\tau},\psi_{\sigma\tau})$ is a simple complex of groups
over $(X_+,V(\tau)_+)$ where $G_{\tau}$ is the local group along the
stratum $V(\tau)_{+}$. Further,
$\varphi=(\varphi_{\tau})_{\tau\in\Delta}:G(\Delta)\ra T_{2}$ is a
simple morphism injective at the local groups.}
\end{prof}

\begin{prof}{\rm $\widehat{G(\Delta)}$ is by definition the free
product of $G_{\tau}$ with the relations
$\psi_{\sigma\tau}(h)=h~\forall~ h\in G_{\tau}$ whenever $\tau<\sigma$
in $\Delta$. Thus, $\widehat{G(\Delta)}$ is simply the graph product
of the vertex groups $G_{\rho_j}\simeq\bz_2$ over the graph $S_N \cap
\Delta(2)$ where the vertices of the graph correspond to $\Delta(1)$
and the edges correspond to $\Delta(2)$. Therefore,
$\widehat{G(\Delta)}\simeq W(\Delta)$ and (3) follows.}\hfill $\Box$\vspace{-.7pc}
\end{prof}
\end{proof}

Let $G$ be a group for which there exists a simple morphism $\phi:{
G(\Delta)}\rightarrow G$, injective at the local groups. Then, $G \times
X_{+}/\sim := \{(g,x):g\in G, x\in X_{+} :(g,x) \sim (g',
x') \Leftrightarrow x=x'; g\cdot (g')^{-1}\in G_{\tau}\}$,
where $\tau$ is the unique cone such that $x\in O_{\tau}$. Let
$D(\Delta,\varphi)=\sqcup_{\tau\in\cq}G/G_{\tau}$. Then,
$D(\Delta,\varphi)$ is a poset consisting of pairs $(g\cdot
G_{\tau},\tau)$ where $\tau\in\Delta$ and $g\cdot G_{\tau}$ is a coset
of $G_{\tau}$ in $G$ and, $D(\Delta,\varphi)$ has the partial order,
$(g\cdot G_{\sigma},\sigma)<(g'\cdot G_{\tau},\tau)$ if and only if
$\sigma<\tau$ in $\Delta$ and $(g')^{-1}\cdot g\in G_{\sigma}$.

\begin{lem}\label{ST} $X$ is a stratified space over
$D(\Delta,\varphi)$. Furthermore{\rm ,} the $T_2$ action on $X$ is
strata-preserving{\rm ,} with $X_+$ as the strict fundamental domain.
\end{lem}

\begin{proof} By definition, $(T_2\times X_{+}/\sim)$ is a stratified
space over $D(\Delta,\varphi)$ such that, the action of $T_2$ on
$T_2\times X_{+}/\sim$ is strata-preserving where, $t\in T_2$ takes the
stratum $(t',V(\tau)_{+})$ to the stratum $(tt',V(\tau)_{+})$. A strict
fundamental domain for this action is the copy $1\times X_{+}$
corresponding to the identity element $1 \in T_2$. However, by
Proposition~\ref{RT}, there is a canonical $T_2$-equivariant isomorphism
from $(T_2\times X_{+})$/$\sim $ to $X$. Thus $X$ gets a structure of a
stratified space over $D(\Delta,\varphi)$ in such a way that, the action
of $T_2$ on $X$ is strata-preserving and the strict fundamental domain
for this action is $X_{+}\subseteq X$.\hfill $\Box$
\end{proof}

\begin{theor}[\!]\label{FG}$\left.\right.$

\begin{enumerate}
\renewcommand\labelenumi{{\rm \arabic{enumi}.}}
\item Let $D(X_+,\varphi)$ and $ D(X_+,\iota)$ denote the developments
of $X_+$ with respect to $\varphi$ and $\iota$ respectively. Then{\rm ,}
$D(X_+,\varphi)\simeq (T_2\times X_{+}/{\sim}) \simeq X $ and
$D(X_+,\iota)\simeq (W \times X_+/{\sim})=: {\widetilde X}$. There are
strata-preserving actions of $W$ on $ D(X_+,\iota) $ and of $T_2$ on
$D(X_+,\varphi)$ with strict fundamental domain $X_{+}$.

\item $X$ is connected if and only if the primitive vectors along the
edges of $\Delta$ span $N\otimes_{\bz}\bz_2$. In particular{\rm ,} $X$ is
connected whenever the primitive vectors contain a $\bz$ basis for $N$.

\item ${\widetilde X}=W \times X_+/{\sim}$ is the universal cover of $X$
and $\pi_{1}(X)\simeq {\rm ker}(\widehat{\varphi})${\rm ,} where{\rm ,}
$\widehat{\varphi}:W \rightarrow T_2$ is the {\it canonical}
homomorphism induced by $\varphi$.

\item Let ${\it h}: W \rightarrow W_{ab}\simeq \bz_{2}^{d}$ be the
surjective group homomorphism obtained by abelianisation. Associated to
the map ${\it h}${\rm ,} we have a simple morphism $\alpha:
{G(\Delta)}\rightarrow \bz_{2}^{d}$ such that $\widehat{\alpha}= {\it
h}$. Then{\rm ,} $D(X_+,\alpha)\simeq \bz_{2}^{d} \times X_+ / \sim$ is a
covering space over $X$ with deck transformation group $\bz_{2}^{d-n}${\rm ,}
it is a covering space of $X$ and $\pi_{1}( D(X_+,\alpha))=[W,W]$.
\end{enumerate}
\end{theor}

\begin{proof}
To prove this theorem we use Prop.~12.20 of \cite{Brid}.

\setcounter{proo}{0}
\begin{prof}{\rm By Prop.~\ref{Coxeter}, the development
$D(X_+,\varphi)$ of $X_+$, with respect to the simple morphism $\varphi$
from the simple complex of groups $G(\Delta)$ over $X_+$ to $T_2$, is a
stratified space over $D(\Delta,\varphi)$ and is isomorphic to
$T_2\times X_+/\sim$ in such a way that, the induced action of $T_2$ on
$T_2\times X_+/\sim$ is identical to that in Lemma \ref{ST}. Hence by
Lemma \ref{ST}, $D(X_+,\varphi)$ is isomorphic to $X$ as a stratified
space and further, the isomorphism is equivariant under the
strata-preserving action of $T_2$. Similarly, the development $D(X_+,\iota)$ of
$X_+$, with respect to the canonical simple morphism $\iota$ from
$G(\Delta)$ to $W$ is isomorphic to $(W\times X_+/\sim)$ which is a
stratified space over the poset $D(\Delta,\iota)$ and further, there is
a strata-preserving action of $W$ on $D(X_+,\iota)$ with strict
fundamental domain, $X_+$.}
\end{prof}

\begin{prof}{\rm From Lemma~\ref{SC}, $X_+$ is contractible, in
particular it is connected. Hence, $D(X_+,\varphi)$ is connected if and
only if $\widehat{\varphi}$ is surjective which is equivalent to the
assumption that the image of the primitive edge vectors span
$N\otimes_{\bz}\bz_2$. This certainly happens if a part of the primitive
vectors along $\Delta(1)$ form a $\bz$-basis for $N$.}
\end{prof}

\begin{prof}{\rm From Lemma \ref{SC} it follows that, $X_+$ is simply
connected and the strata of $X_+$ are arcwise connected. Further, from
Prop.~\ref{Coxeter} we know that $\widehat{G(\Delta)}\simeq W$.
Therefore, $~(W\times X_+/\sim) \simeq D(X_+,\iota) $ is the universal
cover of $X\simeq D(X_+,\varphi)$ and ${\rm ker}(\widehat{\varphi})\simeq
\pi_1(X)$ where, $\widehat{\varphi}$ is the canonical surjective
homomorphism induced by $\varphi$.}
\end{prof}

\begin{prof}{\rm Since $\alpha={\it h}\circ\iota$ and $T_2$ being
abelian $\widehat{\varphi}$ factors through $W_{ab}$. Therefore, the
simple morphism $\alpha$ is injective at the local groups and the
development $D(X_+,\alpha)\simeq\bz_2^d\times X_+/{\sim}$. The
remaining claims of (4) follow simply by the direct application of
Prop.~12.20 of \cite{Brid}.}\hfill $\Box$\vspace{-.7pc}
\end{prof}
\end{proof}

\begin{rem}\label{Conn}{\rm \hskip -.55pc({\em Connectedness of $X$}).\ \ If the primitive edge
vectors of $\Delta(1)$ do not span $N\otimes_{\bz} \bz_2$, then $X$
is not connected and the number of connected components of $X$ is
equal to $[N\otimes_{\bz} \bz_2:\varphi(W)]$.  In fact, $\Delta$ is
supported on a smaller dimensional lattice and therefore, $X$ is
isomorphic to $X^{'}\times (\br^*)^{(k/2)}$, where $k = [N\otimes_{\bz}
\bz_2:\varphi(W)]$ and $X'$ is a connected toric variety of dimension
$n-k$. For example, the real toric variety associated to the fan
$\Delta = \{e_1,-e_1,\{0\}\}$ in $N = \bz e_1\oplus\bz e_2$ is
homeomorphic to $\bs^1\times \br^*$ and has two connected components
$\bs^1\times \br^{+}$ and $\bs^1\times \br^{-}$.  Indeed, for $X$ to
be connected it is not necessary that the primitive edge vectors
should span $N$ for example the real toric variety associated to the
fan $\Delta = \{\langle 2e_1+3e_2\rangle,\langle e_1\rangle,\{0\}\}$ in
$N = \bz e_1\oplus\bz e_2$, is smooth and connected but the edge vectors
$\{2e_1 + 3e_2,e_1\}$ do not form a $\bz$-basis\break for $N$.}
\end{rem}

\section{A presentation for ${\rm p}_{\bf 1}$({\bi X})}

Let $X$ be smooth and connected. In this section we shall give a
presentation for $\pi_1(X)$ with generators and relations defined
purely from the combinatorial structure of $\Delta$.

Let $\{v_{1},\ldots,v_{n}\}$ be primitive vectors along $\Delta(1)$
which form a basis for $N\otimes_{\bz}\bz_{2}$ and let
$\{u_1,\ldots,u_n\}$ be the dual basis. Let $a_{j,i}=\langle
u_i,v_j\rangle~ {\rm mod}~ \bz_{2}$ for $1\leq j\leq d$ and $1\leq i\leq
n$. Then, $A=(a_{j,i})$ is the {\it characteristic matrix of $\Delta$
with respect to $\{v_1,\ldots, v_n\}$}.

For ${\ut} = (t_1,\ldots,t_n) \in \bz_2^n $, let $ b_i^j = t_i +
a_{j,i}$ for $1 \leq i \leq n, 1 \leq j \leq d $ and let
$c^{p,q}_{i}=t_{i}+a_{p,i}+a_{q,i}$ for $ 1\leq i\leq n;~1\leq p,q\leq
d$. We shall denote the vector $(b^j_i)_{i=1,\ldots,n}$ by $\ub^j$ and
the vector $(c^{p,q}_i)_{i=1,\ldots,n}$ by $\uc^{p,q}$.

In the following proposition we will give a presentation for $\pi_1(X)$
using the above\break data.

\setcounter{theore}{0}
\begin{propo}$\left.\right.$\vspace{.5pc}\label{pres}

\noindent The fundamental group $\pi_1(X)$ has a presentation with
generators
\begin{align*}
\{y_{j,\ut} : 1\leq j\leq d\mid \ut=(t_1, \ldots, t_n)\in \bz_2^n\}
\end{align*}
and relations
\begin{align*}
&\bigcup_{{\ut} \in \bz_{2}^n}\{y_{1, (0, \ldots, 0)}^{t_1} \cdot
y_{2,(t_1,0,\ldots,0)}^ {t_2}\cdots y_{n, (t_1, \ldots, t_{n-1},
0)}^{t_n}\}\\[.2pc]
&\bigcup_{\ut\in\bz_2^n}\{ y_{1,{\ut}}\cdot y_{j,\ub^j}\mid 1\leq j\leq
d \} \\[.2pc]
&\bigcup_{\ut\in\bz_2^n}\{ y_{p,{\ut}}\cdot y_{q,\ub^p}\cdot
y_{p,\uc^{p,q}}\cdot y_{q,\ub^q}\mid \langle v_p,v_q\rangle \in\Delta\}.
\end{align*}
\end{propo}

\begin{proof} We know from Theorem~\ref{FG} that $\pi_{1}(X)$ is
isomorphic to the kernel of the surjective homomorphism
${\widehat{\varphi}}:W \rightarrow T_2 \simeq \bz_{2}^{n}$, where $W$
has the presentation $\langle S\mid R\rangle$ for
$S=\{s_1,s_2,\ldots,s_d\}$ and
$R=\{s_1^2,s_2^2,\ldots,s_d^2;~(s_is_j)^2$ whenever $\{v_i,v_j\}$
spans a cone in $\Delta\}$.

We further have the following commuting diagram:
\begin{equation*}
\begin{array}{ccccccccc}
1 & \rightarrow & F^{'} & \rightarrow & F(S) & \rightarrow & \bz_2^n &
\rightarrow & 1\\
& &\downarrow & & \downarrow \psi & &\| & &\\
1 & \rightarrow & H & \rightarrow & W & \rightarrow & \bz_2^n
&\rightarrow & 1
\end{array},
\end{equation*}
where $F(S)$ denotes the free group on $S, \psi$ denotes the canonical
surjection from $F(S)$ to $W, H$ denotes $\pi_{1}(X)$ and $F^{'}
=\psi^{-1}(H)$.

Since ${\mathcal T}=\{s_1^{t_1}\cdot s_2^{t_2}\cdots s_n^{t_n} \mid
(t_1,t_2,\ldots,t_n)~\in~\bz_2^n\}$ is a {\it Schreier transversal} for
$F'$ in $F(S)$, we can apply the Reidemeister--Schreier theorem (cf.
\cite{Cohen,LS}) to obtain a presentation for $\pi_1(X)$ from
that of $W$. Let
\begin{align}
S_H &= \{y_{j,\ut} :1\leq j\leq d; \ut\in \bz_2^n \},\\[.2pc]
R^{1}_H &= \{\alpha_{0}(u)~ \forall~ u\in {\mathcal T}\},\\[.2pc]
R_{\ut} &= \{\alpha_{\ut}(r)~ \forall ~r\in R~;~ {\ut}\in \bz_2^n\},
\end{align}
where $0:=(0,0,\ldots, 0) \in {\bz}_{2}^n$, $\alpha_{\ut}$ :
$F(S)\rightarrow F(S_{H})$ is defined recursively as follows:

\noindent $\alpha_{\ut}(1):=1$; $\alpha_{\ut}(s_j)=y_{j, \ut}$. Suppose
that by induction we have defined $\alpha_{\ut}(w)$ for $w\in F(S)$
then, $\alpha_{\ut}(w\cdot s_j):= \alpha_{\ut}(w)\cdot
\alpha_{\underline{t\cdot s_j}}(s_j)$ where, $\underline{t\cdot s_j}\in
\bz_2^n$ corresponds to the coset representative
$\varphi(w')\in{\mathcal T}$ of $F'\cdot w'$, where $w'=s_1^{t_1}\cdots
s_n^{t_n}\cdot s_j$.

Note that, $~\forall~{\ut} = (t_1,t_2,\ldots,t_n) \in \bz_2^n $ we have
\begin{enumerate}
\renewcommand\labelenumi{(\roman{enumi})}
\leftskip .3pc
\item $\alpha_{0}(s_1^{t_1} \cdot s_2^{t_{2}} \cdots s_n^{t_n})
~=~(\alpha_{0}(s_1))^{t_1} \cdot(\alpha_{(t_1,0,0,\ldots,0)}(s_2))^{t_2}
\cdots(\alpha_{(t_1,t_2,\ldots, t_{n-1},0)}(s_n))^{t_n}$,\vspace{.2pc}

\item $\alpha_{\ut}(s_j^2)~=~\alpha_{\ut}(s_j)\cdot
\alpha_{\ub^j}(s_j)$,\vspace{.2pc}

\item $\alpha_{\ut}((s_p\cdot s_q)^2)~=~\alpha_{\ut}(s_p)\cdot
\alpha_{\ub^p}(s_q)\cdot \alpha_{\uc^{p,q}}(s_p)\cdot
\alpha_{\ub^q}(s_q) ~\forall~1\leq j\leq d~.$
\end{enumerate}
It follows from definition~(3.2) and from the identity (i) above
that
\begin{align*}
R^{1}_H &= \left\{ \begin{array}{@{}ll@{}}
\alpha_{0}(s_1^{t_1}\cdot s_2^{t_2}\cdots s_n^{t_n}) &~\mid~{\ut} =
(t_1,t_2,\ldots, t_n) \in \bz_2^n
\end{array} \right\}\\[.2pc]
&= \left\{ \begin{array}{@{}ll@{}}
(y_{1,(0,0,\ldots,0)}^{t_1} \cdot y_{2,(t_1,0,\ldots,0)}^{t_2}\ldots
y_{n,(t_1,t_2,\ldots, t_{n-1},0)}^{t_n}) \\[.4pc]
|{\ut} = (t_1,t_2,\ldots,t_n) \in \bz_2^n
\end{array}\!\!\right\}.
\end{align*}
Also the definition (3.3) and the identities (ii) and (iii) above, imply
that
\begin{align*}
R_{\ut} &= \left\{ \begin{array}{@{}llll@{}}
\alpha_{\ut}(s_1^2), \ldots,
\alpha_{\ut}(s_d^2)~;\\[.4pc]
(\alpha_{\ut}(s_ps_q)^2) ~{\rm whenever}~\{v_p,v_q\}~{\rm
spans~a~cone~in}~\Delta
\end{array}\!\!\right\}\\[.2pc]
&= \left\{ \begin{array}{@{}llll@{}}
y_{1,{\ut}}\cdot y_{1,\ub^1}, \ldots, y_{d,{\ut}}\cdot
y_{d,\ub^d}~;\\[.4pc]
y_{p,{\ut}}\cdot y_{q,\ub^p)}\cdot y_{p,\uc^{p,q}}\cdot
y_{q,\ub^q}\\[.4pc]
{\rm whenever}~ \{v_p,v_q\}~ {\rm spans~a~cone~in}~ \Delta
\end{array}\!\!\right\}.
\end{align*}
Here, $b^j_i~=~t_i+a_{ji} ~\forall~ 1\leq j \leq d;~1\leq i\leq n$ and
$c^{p,q}_i~=~t_i+a_{pi}+a_{qi}~\forall ~1\leq i\leq n;~1\leq p,q\leq d$.

Let $ R^{2}_H:=~\cup_{\ut \in \bz_2^n} R_{\ut}$. Then, from the
Reidemeister--Schreier theorem it follows that, $\pi_1(X)$ has the
presentation $\langle S_H \mid R^1_H;~ R^2_H \rangle $.\hfill $\Box$
\end{proof}

\begin{lem}\label{spres} The presentation can be simplified with lesser
number of generators and relations by expressing the $y_{j,\ut} \in
S_{H}$ as words in $S$ and using the relations in $W$ if we make further
assumption on $\Delta$ that{\rm ,} $\langle v_p,v_q\rangle \in\Delta$
for all $1\leq p,q\leq n$.
\end{lem}

\begin{proof} $\alpha_{\ut}(s_j)=s_1^{t_1}\cdots s_n^{t_n}\cdot s_j
\cdot \varphi(s_1^{t_1}\cdots s_n^{t_n}\cdot s_j)$ where
$\varphi(w)\in{\mathcal T}$ is the coset representative of $F'\cdot w$
for $w\in W$. Notice further that in $W$, the $\{s_1,\ldots, s_n\}$
commute among themselves by our assumption on $\{v_1,v_2,\ldots,v_n\}$.
Hence we have, $\varphi(s_1^{t_1}\cdots s_n^{t_n}\cdot
s_j)=s_1^{t_1+a_{j,1}}\cdot s_2^{t_2+a_{j,2}}\cdots s_n^{t_n+a_{j,n}}$,
$~\forall~(t_1,t_2,\ldots,t_n)\in~\bz_2^n$. This implies that we have
the following identities:
\begin{enumerate}
\renewcommand\labelenumi{{\rm (\roman{enumi})}}
\leftskip .3pc
\item $y_{j,\ut} ~=~ s_1^{t_1} \cdot s_2^{t_2} \cdots s_n^{t_n}
\cdot (s_j) \cdot s_1^{t_1} \cdot s_2^{t_2} \cdots s_j^{t_j+1} \cdots
s_n^{t_n}~ = 1$ for all $1\leq j\leq n; \ut\in \bz_2^n$ (since,~
$a_{j,i}~=~\delta_{j,i}$~ for all~ $1\leq j\leq n$),\vspace{.2pc}

\item $y_{j,0}~=~\alpha_{0}(s_j)~=~s_j\cdot s_1^{a_{j,1}}\cdot
s_2^{a_{j,2}}\cdots s_n^{a_{j,n}}$~ for all~ $n+1\leq j\leq d$,\vspace{.2pc}

\item $ y_{j,\ut}=~\alpha_{\ut}(s_j)~=~s_1^{t_1}\cdot
s_2^{t_2}\cdots s_n^{t_n}\cdot (s_j\cdot s_1^{a_{j,1}}\cdot
s_2^{a_{j,2}}\cdots s_n^{a_{j,n}})\cdot s_1^{t_1}\cdot
s_2^{t_2}\cdots s_n^{t_n}$ for all $n+1\leq j\leq d$,\vspace{.2pc}

\item $y_{j,\ub^j}=\alpha_{\ub^j}(s_j)= s_1^{t_1+a_{j,1}}\cdots
s_n^{t_n+a_{j,n}}\cdot (s_j)\cdot s_1^{t_1}\cdots s_n^{t_n} =
(s_1^{t_1}\cdots s_n^{t_n}\cdot(s_j\cdot s_1^{a_{j,1}}\cdots
s_n^{a_{j,n}})\cdot s_1^{t_1}\cdots s_n^{t_n})^{-1}=
(y_{j,{\ut}})^{-1}$.
\end{enumerate}
Let $S_j~=~s_j\cdot s_1^{a_{j,1}}\cdots s_n^{a_{j,n}}$ then, from (ii)
and (iii) above, we see that the generators of $\pi_1(X)$ are,
$\{s_1^{t_1}\cdot s_2^{t_2}\cdots s_n^{t_n} \cdot S_j\cdot
s_1^{t_1}\cdot s_2^{t_2}\cdots s_n^{t_n}~ \forall ~n+1\leq j\leq d;
~\forall~ \ut \in \bz_2^n \}$.  Further, since $R^1_H$ consists of
words in $\{y_{j,\ut}$ for $1\leq j\leq n$ and $\ut\in\bz_2^n\}$, from
(i) we see that $R^1_H=\{1\}$. Furthermore, (iv) implies that the
first set of relations in $R_{\ut}$ are of the form $\{y_{j,k}\cdot
(y_{j,k})^{-1}\}$, therefore they trivially hold in the group
$\pi_1(X)$.  Thus, finally, the number of generators reduce to
$(d-n)\cdot 2^n$ and they are: $\{y_{j,\ut}$ for $n+1\leq j\leq
d;~\ut\in\bz_2^n\}$ and the non-trivial relations are of the form:
$R_{\ut}=\{y_{p\varphi({\ut})}\cdot y_{q\varphi(\ub^p)}\cdot
y_{p\varphi(\uc^{p,q})}\cdot y_{q\varphi(\ub^q)}$ whenever
$\{v_p,v_q\}$ spans a cone in $\Delta\} ~\forall~\ut \in
~\bz_{2}^n$. Therefore, the final presentation for $\pi_1(X)$ is
$\langle S_{H},R_{H}\rangle$, where
\begin{align*}
S_{H} &=\{y_{j,\ut}~{\rm where}~n+1\leq j\leq d~{\rm and}~\ut\in\bz_2^n\},\\[.2pc]
R_{\ut}&=\{y_{p,\varphi({\ut})}\cdot y_{q,\varphi(\ub^p)}\cdot
y_{p,\varphi(\uc^{p,q})}\cdot y_{q,\varphi(\ub^q)}\\[.2pc]
&\quad ~{\rm whenever}~ \{v_p,v_q\}~{\rm spans~ a~ cone~ in}~ \Delta\}
~\forall~\ut \in \bz_{2}^n,\\[.2pc]
R_{H} &= \cup_{\ut\in\bz_2^n}R_{\ut}.
\end{align*}
Further, $\pi_1(X)$ is generated as a subgroup of $W(\Delta)$ by
$s_1^{t_1}\cdots s_n^{t_n}\cdot S_j\cdot s_1^{t_1}\cdots s_n^{t_n}$,
where $\ut\in \bz_2^n$ and $S_j=s_j\cdot s_1^{a_{j,1}}\cdots
s_n^{a_{j,n}}$ for $n+1\leq j\leq d$.\hfill $\Box$
\end{proof}

\begin{rem}{\rm By the classification of two-dimensional smooth complete fans
(cf. p.~42 of \cite{Fult} and p.~57 of \cite{Jurk}) we observe that,
except the torus $\bs^1\times \bs^1$, all other smooth complete real
toric surfaces correspond bijectively to the two-dimensional compact
non-orientable manifolds. This is because, they are obtained by
successively blowing up $\bp_{\br}^2$ at $T$-fixed points and are
therefore homeomorphic to $\bp_{\br}^2\#\cdots\#\bp_{\br}^2$ (${d-2}$
copies). However, the classical presentation for the fundamental group
is apparently different from the presentation we have obtained,
especially because it has only one relation. In the cases when $d=3$
and 4, where the spaces are $\bp_{\br}^2$, $\bs^1\times \bs^1$ or
the Klein-bottle$~\simeq \bp_{\br}^2\#\bp_{\br}^2$, the presentations
we give agrees with some of the classical ones.  We hope to simplify
the presentation given above to reduce the number of generators and
relations so that in general it agrees with the classical cases.}
\end{rem}

\begin{rem}{\rm Note that the fundamental group $\pi_1(X)$ and hence its
presentation depends only on the $2$-skeleton $\Delta(2)$ of $\Delta$.}
\end{rem}

\section{The coxeter group {\bi W}($\Delta$)}

In this section we prove some general results on right-angled Coxeter
groups and in particular for $W(\Delta)$. Let $M=(m_{ij})$ denote the
Coxeter matrix corresponding to $W$.

\setcounter{theore}{0}
\begin{lem}\label{ab}$[W,W]$ is abelian if and only if for all $1\leq
j\leq d$ there exists at most one $i$ such that $\langle v_i,v_j\rangle
\notin \Delta$.
\end{lem}

\begin{proof} If there exist $i\neq k$ such that $\{v_i,v_j\}$ and
$\{v_k,v_j\}$ does not span a cone in $\Delta$ then, $[s_i,s_j]\cdot
[s_k,s_j]\neq [s_k,s_j]\cdot[s_i,s_j]$ in $[W,W]$.

Conversely, if for each $1\leq j\leq d$, there exists at most one $i$
such that $\langle v_i,v_j\rangle\notin \Delta$ then, using the
relations in $W$, it is easy to see that for any word $w\in W$,
$w\cdot [s_i,s_j]\cdot w^{-1}=[s_i,s_j]$ or $[s_j,s_i]$. (It is
$[s_j,s_i]$ iff either one of $s_i$ or $s_j$ but not both occurs in
the reduced expression of $w$.) Now, $[W,W]$ is the normal subgroup of
$W$ generated by the commutators $\{[s_i,s_j]\mid \langle
v_i,v_j\rangle\notin\Delta\}$.  Therefore, under the above assumption,
$\{[s_i,s_j]\mid \langle v_i,v_j\rangle \notin \Delta\}$ in fact
generate $[W,W]$ as a subgroup of $W$. Further, since they commute
among themselves, $[W,W]$ is abelian.\hfill $\Box$
\end{proof}

\begin{lem}\label{FO} A word $w\in W$ is of finite order if and only if
it is of order $2$. Moreover{\rm ,} in this case{\rm ,} $w$ is a conjugate in $W$ to
a word $w'$ which is of the form{\rm ,} $w'=s_{j_1}\cdots s_{j_l}$ with
$s_{j_p}\cdot s_{j_q}=s_{j_q}\cdot s_{j_p}~\forall~1\leq p,q\leq l$.
\end{lem}

\begin{proof} Suppose $w=v\cdot w'\cdot v^{-1}$ where, $w'$ is as above
and $v\in W$. Then $w$ is clearly of order $2$. On the other hand, if
$w$ is not of the above form then the reduced expression for $w$ is of
the form $w=s_{i_1}\cdots s_{i_k}$ where $s_{i_p}$ and $s_{i_q}$ do not
commute for some $1\leq p,q\leq k$. Indeed, by repeatedly using the
relation $s_i\cdot s_j= s_j\cdot s_i$ whenever $m_{ij}=2$, we can assume
without loss of generality that, up to conjugation $w$ is of the form
$s_{i_1}\cdot s_{i_2}\cdots s_{i_k}$, where $s_{i_1}$ and $s_{i_k}$ do
not commute. Then it follows that, for any positive integer $r$,
$w^r=(s_{i_1}\cdots s_{i_k})\cdot(s_{i_1}\cdots
s_{i_k})\cdots(s_{i_1}\cdots s_{i_k})$ is in fact a reduced expression
in $W$. Hence, $w$ is of infinite order.\hfill $\Box$
\end{proof}

\begin{lem}\label{flag} Let $w=s_{i_1}\cdots s_{i_k}\in W${\rm ,} where
$\langle v_{i_1},\ldots,v_{i_k}\rangle\in\Delta $ and let
$w'=s_{j_1}\cdots s_{j_l}$ where{\rm ,} $\langle
v_{j_p},v_{j_q}\rangle\in\Delta$ for all $1\leq p,q\leq l$ but $\langle
v_{j_1},\ldots, v_{j_l}\rangle\notin\Delta$. Then{\rm ,} $w\notin N(w')${\rm ,}
where $N(w')$ is the normal subgroup generated by $w'$ in $W$.
\end{lem}

\begin{proof}
Suppose on the contrary $w=v_1\cdot w'\cdot v_1^{-1}\cdot v_2\cdot
w'\cdot v_2^{-1}\cdots v_r\cdot w'\cdot v_r^{-1}$ for some
$v_1,v_2,\ldots v_r\in W$. By Lemma~\ref{FO} we know that $(w')^{2}=1$.
Hence, the above expression can be rewritten as
\begin{enumerate}
\renewcommand\labelenumi{\arabic{enumi}.}
\item $w=[v_1,w']\cdot[w',v_2]\cdot[v_3,w']\cdots[w',v_r]$ if $r$ is
even,

\item $w=[v_1,w']\cdot[w',v_2]\cdots [v_r,w']\cdot w'$ if $r$ is odd.\vspace{-.4pc}
\end{enumerate}
This implies that, $w\in[W,W]$ if $r$ is even and $w\cdot w'=w\cdot
(w')^{-1}\in [W,W]$ if $r$ is odd.

Now let $h:W\rightarrow\bz_2^d$ be the abelianisation map which takes
$s_j$ to the coordinate vector $e_j=(0,0,\ldots,1,\ldots,0)$ (with 1 at
the $j$th position). Also, by our choice of $w$ and $w'$ we observe
that, $\{s_{i_1},\ldots,s_{i_k}\}$ and $\{s_{j_1},\ldots,s_{j_l}\}$
pairwise commute in $W$ and the tuples $(i_1,\ldots, i_k)$ and
$(j_1,\ldots,j_l)$ are distinct.

Therefore, $h(w)=\Sigma_{p=1}^{k} e_{i_p}\neq (0,\ldots,0)$ when $r$ is
even and $h(w\cdot w')=\Sigma_{p=1}^{k}e_{i_p} +
\Sigma_{q=1}^{l}e_{j_q}\neq(0,\ldots,0) $ when $r$ is odd. This is a
contradiction since on the other hand, $w$ and $w\cdot (w^{'})^{-1}\in
[W,W]$ when $r$ is even and $r$ is odd respectively. This proves the
lemma.

\hfill $\Box$
\end{proof}

\begin{rem}{\rm Lemma~\ref{ab} if phrased differently will be: $[W,W]$
is abelian if and only if there exists at most one $i$ for every $j$
such that $m_{i,j}\neq 2$, holds not just for right-angled Coxeter
groups, but for more general class of Coxeter groups with
$m_{i,j}=2$ or $m_{i,j}\geq 5~\forall~ i,j$.}
\end{rem}

\section{Criterion for ${\rm p}_{\bf 1}$({\bi X}) to be abelian}

Let $X$ be smooth and connected. In the following theorem we give
conditions on $\Delta$ under which $\pi_1(X)$ is abelian. We shall
follow the notations in \S3 and further assume that $\langle
v_p,v_q\rangle\in\Delta$ for every $1\leq p,q\leq n$ as in
Lemma~\ref{spres}.

\setcounter{theore}{0}
\begin{theor}[\!]\label{abel} $\pi_1(X)$ is abelian if and only if one
of the following holds in $\Delta$.
\begin{enumerate}
\renewcommand\labelenumi{{\rm \arabic{enumi}.}}
\item For every $1\leq i,j\leq d$, $\{v_i,v_j\}$ spans a cone in
$\Delta$. In this case{\rm ,} $\pi_1(X)$ is isomorphic to $\bz_2^{d-n}$.

\item For each $~1\leq j\leq d$ there exists at most one $i=i_j$ with
$1\leq i_j\leq n$ such that{\rm ,} $\{v_{i_j},v_j\}$ does not span a cone in
$\Delta$ and $\langle u_{i_j},v_j\rangle =1~{\rm mod}~2$. Further{\rm ,} for each
$n+1\leq k\leq d~$ such that $k \neq j$ we have{\rm ,} $\langle
u_{i_j},v_k\rangle=0~{\rm mod}~2$.
\end{enumerate}
\end{theor}

\begin{proof} Recall that we have an exact sequence, $1\ra
[W,W]\ra \pi_1(X)\ra\bz_2^{d-n}\ra 1$ and further, $[W,W]$ is
generated as a normal subgroup of $W$ by $[s_{i_p},s_{i_q}]$ whenever
$\{v_{i_p},v_{i_q}\}$ does not span a cone in $\Delta$.

\begin{step}{\rm Since $[W,W]$ is a subgroup of $\pi_1(X)$, if $\pi_1(X)$
is abelian then $[W,W]$ must be abelian. By Lemma \ref{ab}, $[W,W]$ is
abelian if and only if for every $v_j$ there exists at most one $v_i$
such that $\{v_j,v_i\}$ does not span a cone in $\Delta$.}
\end{step}

Further, $[W,W]=\{1\}$ if and only if any two $\{v_j,v_i\}$ for $1\leq
i,j\leq d$ spans a cone in $\Delta$ which implies that, $W\simeq\bz_2^d$
and $\pi_1(X)\simeq\bz_2^{d-n}$.

\begin{step}{\rm On the other hand, if $[W,W]\neq \{1\}$ then there
exists $\{v_j,v_i\}$ which does not span a cone in $\Delta$. However,
since $[W,W]$ is abelian, this $i=i_j$ must be unique for every
such $j$. Thus, in $W$, $s_j$ and $s_{i_j}$ do not commute but they
both commute with $s_k$ for every $1\leq k\leq n$.}
\end{step}

\begin{step}{\rm Suppose now that, for some $n+1\leq j\leq d$
we have $n+1\leq i_j\leq d$, then $\pi_1(X)$ is non-abelian for if
$S_j$ denotes the word $s_j\cdot s_1^{a_{j,1}}\cdots s_n^{a_{j,n}}$ in
$W$ then,
\begin{align*}
S_j\cdot S_{i_j} &= s_j\cdot s_{i_j}\cdot s_1^{a_{j,1}}\cdots
s_n^{a_{j,n}}\cdot s_1^{a_{i_j,1}}\cdots s_n^{a_{i_j,n}}\\[.2pc]
&\neq s_{i_j}\cdot s_j\cdot s_1^{a_{j,1}}\cdots s_n^{a_{j,n}}\cdot
s_1^{a_{i_j,1}}\cdots s_n^{a_{i_j,n}} = S_{i_j}\cdot S_j.
\end{align*}

Hence if $\pi_1(X)$ is abelian then, for every $n+1\leq j\leq d$ there
is a unique index $i_j$ such that $\langle v_j,v_{i_j}\rangle\notin
\Delta$ and further, $1\leq i_j\leq n$.}
\end{step}

\begin{step}{\rm Now, if for some $n+1\leq k\leq d$ with $k\neq j$ we
have $a_{k,i_j}=\langle u_i,v_{i_j}\rangle~{\rm mod}~ \bz_2=1$, then
$\pi_1(X)$ is non-abelian. This is because, if $w=[S_k,S_j]\in
\pi_1(X)$ then $w\neq 1$ which we can see by the following cases:
\begin{enumerate}
\renewcommand\labelenumi{}
\item If $a_{k,i_k}=0$ and $a_{j,i_k}=0$, then $w=[s_{i_j},s_j]\neq 1$.\vspace{.2pc}

\item If $a_{k,i_k}=1$ and $a_{j,i_k}=1$, then $w=[s_{i_j},s_j]\neq 1$.\vspace{.2pc}

\item If $a_{k,i_k}=1$ and $a_{j,i_k}=0$, then $w=[s_{i_j},s_j]\cdot
[s_k,s_{i_k}]\neq 1$.\vspace{.2pc}

\item If $a_{k,i_k}=0$ and $a_{j,i_k}=1$, then $w=[s_{i_j},s_j]\cdot
[s_k,s_{i_k}]\neq 1$.
\end{enumerate}

(Here we omit the proofs of the assertion that $w\neq 1$ in each case,
as it follows easily from the relations in $W$).}
\end{step}

\begin{step}{\rm If $a_{j,i_j}=0$, then again $\pi_1(X)$ is non-abelian
since, the elements $s_{i_j}\cdot S_j\cdot s_{i_j}$ and $ S_j^{-1}$ do
not commute in $\pi_1(X)$. This is because, by Lemma \ref{FO}
$[s_{i_j},s_j]$ is an element of infinite order in $W$ and hence,
$(s_{i_j}\cdot S_j\cdot s_{i_j})\cdot S_j^{-1}= [s_{i_j},s_j]\neq
[s_j,s_{i_j}]= S_j^{-1}\cdot (s_{i_j}\cdot S_j\cdot s_{i_j})$.}
\end{step}

\begin{step}{\rm Therefore, if $\pi_1(X)$ is abelian and $[W,W]\neq
1$, then it is necessary that the following conditions must hold:

For every $1\leq j\leq d$, there exists a unique index $i_j$ with $1\leq
i_j\leq n$ such that $\{v_j,v_{i_j}\}$ does not span a cone in $\Delta$
and $a_{j,i_j}=1$. Further, for every $n+1\leq k\leq d$ such that $k\neq
j$, we have $a_{k,i_j}=0$.

We shall now prove that these conditions are in fact sufficient for
$\pi_1(X)$ to be abelian.}
\end{step}

\begin{claim}$\left.\right.$

\begin{enumerate}
\renewcommand\labelenumi{(\roman{enumi})}
\leftskip .2pc
\item $S_j$ and $S_k$ commute for $n+1\leq j,k\leq d$.

\item $w\cdot S_j\cdot w^{-1}$ and $S_j$ commute where,
$w=s_1^{t_1}\cdots s_n^{t_n}$ for every $\ut=(t_1,\ldots,t_n)\in
\bz_2^n$ and $n+1\leq j\leq d$.
\end{enumerate}
\end{claim}

\begin{potc}$\left.\right.$

\begin{enumerate}
\renewcommand\labelenumi{(\roman{enumi})}
\leftskip .2pc
\item $\left.\right.$\vspace{-1.45pc}
\begin{align*}
\hskip -.5cm S_j\cdot S_k &= s_j\cdot s_k\cdot s_1^{a_{j,1}}\cdots s_n^{a_{j,n}}\cdot
s_1^{a_{k,1}}\cdots s_n^{a_{k,n}}\{{\rm since}~ a_{j,i_k}=0~ {\rm
by~assumption}\}\\[.2pc]
\hskip -.5cm &= s_k\cdot s_j\cdot s_1^{a_{j,1}}\cdots s_n^{a_{j,n}}\cdot
s_1^{a_{k,1}}\cdots s_n^{a_{k,n}}\{ {\rm since}~ k\neq i_j\}\\[.2pc]
\hskip -.5cm &= s_k\cdot s_1^{a_{k,1}}\cdots s_n^{a_{k,n}}\cdot s_j\cdot
s_1^{a_{j,1}}\cdots s_n^{a_{j,n}}\{{\rm since}~ a_{k,i_j}=0~ {\rm
by~assumption}\}\\[.2pc]
\hskip -.5cm &= S_k\cdot S_j.
\end{align*}

\item Let $w=s_1^{t_1}\cdots s_n^{t_n}$ such that
$(t_1,t_2,\ldots,t_n)\in \bz_2^n$.
\begin{equation*}
\hskip -.5cm w\cdot S_j\cdot w^{-1}=(s_1^{t_1}\cdots s_n^{t_n})\cdot
(s_j\cdot s_1^{a_{j,1}}\cdots s_n^{a_{j,n}})\cdot(s_1^{t_1}\cdots
s_n^{t_n})= z.
\end{equation*}
\end{enumerate}
\begin{enumerate}
\renewcommand\labelenumi{(\alph{enumi})}
\item If $t_{i_j}=0$, then $z=s_j\cdot s_1^{a_{j,1}}\cdots
s_n^{a_{j,n}}=S_j$ \{since~$s_j\cdot w=w\cdot s_j$\}. Thus $w\cdot
S_j\cdot w^{-1}=S_j$.

\item If $t_{i_j}=1$, then $z=s_{i_j}\cdot s_j\cdot s_1^{a_{j,1}}\cdots
s_n^{a_{j,n}}\cdot s_{i_j}~=s_{i_j}\cdot S_j\cdot s_{i_j}$. Further,
\begin{align*}
\hskip -.5cm (w\cdot S_j\cdot w^{-1})\cdot (S_j) &= s_{i_j}\cdot S_j\cdot s_{i_j}\cdot
S_j\\[.2pc]
\hskip -.5cm &= s_{i_j}\cdot s_j\cdot s_1^{a_{j,1}}\cdots \widehat{s_{i_j}}\cdots
s_n^{a_{j,n}}\cdot S_j\ \{{\rm since}~~a_{j,i_j}=1\}\\[.2pc]
\hskip -.5cm &= s_{i_j}\cdot s_j^2\cdot s_{i_j}=1.
\end{align*}
This implies that, $w\cdot S_j\cdot w^{-1}~=~S_j^{-1}$.
\end{enumerate}
Hence, $w\cdot S_j\cdot w^{-1}$ commutes with $S_k$ for all $~n+1\leq
j, k\leq d$, since we have either $w\cdot S_j\cdot w^{-1}=S_j$ or
$S_j^{-1}$ in each of the cases. Therefore, since the generators
commute among themselves, we conclude that $\pi_1(X)$ is
abelian.\hfill $\Box$
\end{potc}
\end{proof}

\begin{rem}{\rm If $\Delta$ is complete, then the condition
$a_{j,i_j}=1$ will be forced after Step~4 in which case, we shall
skip Step~5. However this is not true in general for example, in
the non-complete fan $\Delta=\{\{0\},\langle e_1,e_2\rangle,\langle
-2e_1+e_2\rangle\}$ in $N=\bz e_1\oplus\bz e_2$.}
\end{rem}

\begin{rem}\hskip -.5pc\label{fo}{\rm ({\em Torsion elements}).\ \ By Lemma~\ref{FO}, since
$\pi_1(X)$ is a subgroup of $W$, the torsion elements in $\pi_1(X)$ are
always of order $2$. In particular, when $\pi_1(X)$ is abelian,
$S_j~=~s_j\cdot s_1^{a_{j,1}}\cdots s_n^{a_{j,n}}$ for $n+1\leq j\leq d$
is of order $2$ iff $\langle v_j,v_i\rangle\in\Delta$ for all $1\leq
i\leq n$ and it is of infinite order iff there exists a unique $1\leq
i_j\leq n$ such that $\langle v_j,v_{i_j}\rangle \notin \Delta$ since in
this case, $a_{j,{i_j}}=1$ and $S_j^{2}=[s_j,s_{i_j}]\neq 1$ in
$[W,W]\subset W$.}
\end{rem}

\begin{rem}\label{abelian}{\rm If $\pi_1(X)$ is abelian then $\pi_1(X)$ is
generated by $S_j~=~s_j\cdot s_1^{a_{j,1}}\cdots s_n^{a_{j,n}}$ for
$n+1\leq j\leq d$. Let $\{j_1,j_2,\ldots,j_r\}=J=\{j\mid n+1\leq j\leq
d$ and $\langle v_j,v_{i_j}\rangle\notin \Delta$ for some $1\leq i_j\leq n
\}$. Therefore, if $j\notin J$ then, $\langle v_j,v_i\rangle\in \Delta $
for every $1\leq i\leq n$. Thus, $\pi_1(X)\simeq \bz_2^{d-n-r}\oplus
\bz^r$ where, $\bz^r = \langle S_{j_p}=s_{j_p}\cdot
s_1^{a_{j_p,1}}\cdots s_n^{a_{j_p,n}}$ for $1\leq p\leq r \rangle$.
Furthermore, $[W,W]=\langle [s_{j_p},s_{i(j_p)}]=S_{j_p}^{2}$ for $1\leq
p\leq r\rangle\subset W$ is a free abelian of rank~$r$. We therefore have
the following commuting diagram
\begin{equation*}
\begin{array}{ccccccccc}
1 &\lr &[W,W] &\lr &\pi_1(X) &\lr &\bz_2^{d-n} &\lr &1\\
  & &\| & & \| & &\| & &\\
1 &\lr &\bz^r &\stackrel{\times 2}{\lr} &\bz^r\oplus \bz_2^{d-n-r} &\lr
&\bz_2^{r}\oplus \bz_2^{d-n-r} &\lr &1
\end{array}.
\end{equation*}}
\end{rem}

\begin{rem}{\rm If $\pi_1(X)$ is abelian then necessarily $d\leq 2n$
because, to every $n+1\leq j\leq d$, we associate a unique $i_j$ with
$1\leq i_j\leq n$. Examples of toric varieties with abelian fundamental
group are: (i) products of real projective spaces, (ii) toric bundles
with base as aspherical toric variety with abelian fundamental group and
fibre $\bp^n_{\br}$ for $n\geq2$. (However, this is not true for a
non-trivial bundle with fibre $\bp^1_{\br}$ for example, $\pi_1((\mathbb
{F}_1)_{\br})$ is non-abelian where $(\mathbb{F}_1)_{\br}$ denotes the
real part of the Hirzebruch surface $\mathbb{F}_1.)$}
\end{rem}

\section{Asphericity of {\bi X}}

Let $S_{N}:=(N_{\br}-\{0\})/\br_{>0}$ and $\pi:N_\br-\{0\}\lr S_{N}$ be
the projection. Let $\cs_{\Delta}$ denote the simplicial complex
associated to the smooth fan $\Delta$, where each $k$-dimensional
$\sigma\in \Delta$ corresponds to a $(k-1)$-dimensional spherical
simplex $\pi(\sigma-\{0\})$. If further we assume $\Delta$ to be
complete, then it gives rise to a triangulation of $S_{N}$ (cf. p.~52 of
\cite{Oda}).

Recall that a simplicial complex ${\mathcal S}$ with vertices ${\mathcal
V}=\{v_{i}\}$ is called a {\it flag complex} if the following condition
holds for every finite subset $\{v_{1},v_{2},\ldots,v_{n}\}$ of ${\mathcal
V}$: If $\{v_{i},v_{j}\}$ span a simplex in ${\mathcal S}$ for all $i,j
\in \{1,2,\ldots,n\}$ then $\{v_{1},v_{2},\ldots,v_{n}\}$ span a simplex of
${\mathcal S}$.

Hence, ${\mathcal S}_{\Delta}$ is a flag complex if and only if for
every collection of primitive edge vectors
$\{v_{i_{1}},\ldots,v_{i_{r}}\}$, if $\{\langle
v_{i_{k}},v_{i_{l}}\rangle \in\Delta~\forall~1\leq k, l \leq r\}$ then
$\langle v_{i_{1}},\ldots,v_{i_{r}}\rangle\in\Delta$. We shall say that
$\Delta$ is {\em flag-like} whenever $\cs_{\Delta}$ is a flag complex.

\setcounter{theore}{0}
\begin{theor}[\!]\label{aspher1}$X$ is aspherical if and only if
$\Delta$ is flag-like.
\end{theor}

\begin{proof} If $\widetilde{X}$ is contractible then we claim that
$\Delta$ is flag-like.

Suppose on the contrary, $\Delta$ is not flag-like. Then, $\exists$
$\{v_{j_1},\ldots,v_{j_l}\}$ such that $\forall~ 1\leq p,q\leq l$,
$\langle v_{j_p},v_{j_q}\rangle\in \Delta$ but, $\langle
v_{j_1},\ldots,v_{j_l}\rangle\notin \Delta$.

Let $w'=s_{j_1}\ldots s_{j_l}\in W$ and let $N(w')$ be the normal
subgroup of $W$ generated by $w'$ as in Lemma~\ref{flag}. Also let
$\theta:W\ra W/N(w')$ be the canonical surjection. Clearly,
$\lambda=(\lambda_{\tau}=\theta\circ \iota_{\tau})$ is a simple morphism
from $G(\Delta)\ra W/N(w')$. Further, Lemma~\ref{flag} implies that,
$\lambda_{\tau}:G_{\tau} \subseteq T_2 \ra W/N(w')$ is injective
$\forall~\tau\in\Delta$. Hence, $\lambda$ is injective at the local
groups. Now, the development $D(X_+,\lambda)$ of $X_+$ with respect to
$\lambda$, has $D(X_+,\iota)\simeq\widetilde{X}$ as the universal cover
and its fundamental group $\pi_1(D(X_+,\lambda)\simeq N(w')$ has $w'$ as
a torsion element. This is a contradiction since $D(X_+,\lambda)$ is a
$K(\pi,1)$ space, because of our assumption that $\widetilde{X}$ is
contractible.

For proving the converse, we apply Corollary~10.3 of the main result of
\cite{Davis} to the reflection system $(\Gamma=W$,$V=S)$ on
$M=\widetilde {X}$ with fundamental chamber $Q=X_+$ (which is
contractible by Lemma \ref{SC}). Here, for every $T\subseteq S$,
$Q_T=\cap_{s_{j}\in T}V(\rho_j)_{+}$. Let $W_T$ be the subgroup
generated by $T$ in $W$. Then, the following statements are equivalent:
\begin{enumerate}
\renewcommand\labelenumi{\arabic{enumi}.}
\item $Q_T$ is acyclic for all $T\subseteq S$ with $W_T$ finite.

\item $\Delta$ is  flag-like.
\end{enumerate}

\setcounter{section}{1}
\setcounter{ppro}{1}
\begin{proff}{\rm Let $\rho_{j_1},\ldots,\rho_{j_l}$
be edges such that $\{\rho_{j_p},\rho_{j_q}\!\}$ spans a cone in $\Delta$
for all $1\leq p,q\leq l$. Then, (1) implies that
$Q_T=\cap_{r=1}^{l}V(\rho_{j_r})_+=V(\tau)_+ $ is non-empty since, by
Lemma \ref{FO}, $W_T=\langle s_{j_i}\ldots s_{j_l}\rangle$ is a finite
subgroup of $W$. This implies that,
$\tau=\langle\rho_{j_1},\ldots,\rho_{j_l}\rangle$ is a non-empty cone
in $\Delta$.}
\end{proff}

\setcounter{section}{2}
\setcounter{ppro}{0}
\begin{proff}{\rm Let $T=\{s_{j_1},\ldots,s_{j_l}\}\subseteq S$ be such
that $W_T$ is finite. Then, in particular, $w'=s_{j_1}\cdots s_{j_l}$ is
an element of finite order in $W$. By Lemma~\ref{FO}, the edge vectors
$v_{j_1}\ldots,v_{j_l}$ pairwise span cones in $\Delta$. The assumption
(2) further implies that, $v_{j_1},\ldots,v_{j_l}$ together span a cone
$\tau$ in $\Delta$. Thus, $Q_{T}=\cap_{r=1}^{l} V(\rho_{j_r})_{+} =
V(\tau)_{+}$ is non-empty. Moreover, $V(\tau)$ being a smooth toric
variety, its non-negative part $V(\tau)_+$ is contractible by
Lemma~\ref{SC}, and is hence acyclic if it is non-empty.}
\end{proff}
We therefore conclude from Corollary~10.3 of \cite{Davis} that, if
$\Delta$ is a flag-like then $M=\widetilde{X}$ is contractible.
\hfill $\Box$
\end{proof}

\setcounter{section}{6}
\begin{rem}{\rm In fact, since $(1)\Leftrightarrow (2)$ above, it is
clear that Corollary~10.3 of \cite{Davis} also proves the first
implication of the above theorem. However, in our particular case (where
$W$ is a right-angled Coxeter group), the argument given above is
self-contained and is an application of the `method of development' which is
consistent with the theme of this paper.}
\end{rem}

The following are some corollaries of the above theorem.

\begin{coro}$\left.\right.$\vspace{.5pc}

\noindent If $X$ is aspherical then $V(\tau)$ is aspherical for every
cone $\tau$.
\end{coro}

\begin{proof} This is immediate because, $V(\tau)$ is the toric variety
associated to the fan ${\rm Star}(\tau)$ which by definition (cf. p.~52
of \cite{Fult}) is smooth and flag-like whenever $\Delta$ is smooth and
flag-like. A proof for this is as follows: Let
$\bar{\rho_{i_1}},\ldots,\bar{\rho_{i_k}}$ be edge vectors which
pairwise span cones in ${\rm Star}(\tau)$. Therefore by the definition
of ${\rm Star}(\tau)$, the edges of $\tau$ and $\rho_{i_1},\ldots
\rho_{i_k}$ pairwise span cones in $\Delta$. Since $\Delta$ is
flag-like, this implies that $\gamma=\langle
\tau,\rho_{i_1}\ldots,\rho_{i_k}\rangle$ is a cone in $\Delta$ and
hence,
$\bar{\gamma}=\langle\bar{\rho_{i_1}},\ldots,\bar{\rho_{i_k}}\rangle$ is
a cone in ${\rm Star}(\tau)$. Thus, ${\rm Star}(\tau)$ is
flag-like.

\hfill $\Box$
\end{proof}

\begin{coro}$\left.\right.$\vspace{.5pc}

\noindent Let $X$ be smooth and complete. We can blow up $X$ along a
number of $T$-stable subvarieties to get a smooth complete toric variety
$X'$ which is aspherical.
\end{coro}

\begin{proof} Since $\Delta$ is a smooth and complete fan,
$\cs_{\Delta}$ is a simplicial decomposition of the sphere $S_N$. It is
known that the barycentric subdivision of any simplicial complex is a
flag complex (cf. \cite{Brid}, p.~210). Therefore, if $\Delta'$ is the
refinement of $\Delta$ obtained by taking the cones over the simplices
in the barycentric subdivision of $\cs_{\Delta}$, then $\Delta'$ is a
flag-like fan. It is not difficult to see that $\Delta'$ is also smooth
and complete. Hence, the smooth complete toric variety $X(\Delta')$
which is obtained by blowing up $X$ along certain $T$-stable
subvarieties is aspherical.\hfill $\Box$
\end{proof}

\begin{rem}{\rm However, in some cases we need lesser number of blow ups
to arrive at an aspherical space. For e.g. (i) $\bp_\br ^{2}$ blown up
at a $T$-fixed point is the Hirzebruch surface $({\Bbb F}_1)_{\br}$ (the
Klein-bottle) and $({\Bbb F}_{1})_{\br}$ is aspherical. (ii) $\bp_{\br}^{2}\times \bs^{1}$ needs to
be blown up along a $T$-stable $\bp_{\br}^1$ to get an aspherical space
$(\Bbb{F}_1)_{\br}\times \bs^1$.}
\end{rem}

\section{Subspace arrangement related to $\Delta$}

Throughout this section we assume that $\Delta$ is a smooth and
complete fan.

In this section we define a real subspace arrangement associated to
$\Delta$ whose complement in $\br^d$ is denoted by $\cc_{\Delta}$.
Recall from \cite{Cox} that, $X_{\bc}\simeq X'_{\bc}/(\bc^{*})^{d-n}$
where, $X'_{\bc}$ is the complement of a complex subspace arrangement in
$\bc^d$. By restricting scalars to $\br$ in the above quotient, we show
that $X\simeq \cc_{\Delta}/(\br^{*})^{d-n}$ where $\cc_{\Delta}\simeq
X'_{\br}$. We compute the fundamental group of $\cc_{\Delta}$ and also
give necessary and sufficient conditions for it to be a $K(\pi,1)$
space.

\setcounter{theore}{0}
\begin{definit}$\left.\right.$\vspace{.5pc}

\noindent {\rm A collection ${\mathcal
P}=\{\rho_{i_1},\rho_{i_2},\ldots,\rho_{i_k}\}$ of edges in $\Delta$ is
called a primitive collection if
$\{\rho_{i_1},\rho_{i_2},\ldots,\rho_{i_k}\}$ together does not span a
cone in $\Delta$ but every proper subcollection of ${\mathcal P}$ spans
a cone in $\Delta$. For the primitive collection $\cp$ let ${\mathcal
A(P)}=\left \{ (x_1\ldots, x_d)\in \br^d~|~ x_{i_1}=x_{i_2}=\cdots
=x_{i_k}=0\right\}$.}
\end{definit}

\begin{definit}$\left.\right.$\vspace{.5pc}\label{coordsub}

{\rm \begin{enumerate}
\renewcommand\labelenumi{(\roman{enumi})}
\leftskip .2pc
\item The coordinate subspace arrangement in $\br^d$ corresponding to a
fan $\Delta$ denoted by $\ca_{\Delta}$ is defined as follows:
${\ca}_{\Delta}=\cup_{\cp}{\ca(\cp)}$, where the union is taken over
all primitive collections ${\cp}$ of edges in $\Delta$.

\item Let ${\cc}_{\Delta}$ denote the complement of ${\ca_{\Delta}}$ in
${\br}^d$, i.e. ${\cc}_{\Delta}:= {\br}^d - {\ca_{\Delta}}$.
\end{enumerate}}
\end{definit}

Let $\{{\mathcal P}_1,{\mathcal P}_2,\ldots,{\mathcal P}_r\}$ be the set
of all primitive collections in $\Delta$ consisting of two edges. Let
${\mathcal P}_i~=~\{\rho_{i_p},\rho_{i_q}\}$ where $1\leq i_p,i_q\leq d
~ ~\forall~ 1\leq i\leq r$.

The following lemma generalizes the description of a smooth complete
complex toric variety as given in \cite{Cox} and \cite{Audin} to the
corresponding real and non-negative parts. Although this follows
almost immediately from the complex case, we give a proof for it
since we have not seen the result mentioned anywhere explicitly.

\begin{lem}\label{GQ} The real toric variety $X$ corresponding to a
smooth complete fan $\Delta$ is the geometric quotient of ${\mathcal
C}_{\Delta}$ by the real algebraic torus $(\br^*)^{d-n}$ and we have a
locally trivial principal bundle with total space $\cc_{\Delta}${\rm ,}
base $X$ and structure group $(\br^*)^{d-n}$, i.e.{\rm ,} ${\mathcal
C}_{\Delta}\ra {\mathcal C}_{\Delta}/(\br^*)^{d-n}\simeq X$.
Similarly{\rm ,} $X_{+} \simeq {(\cc_{\Delta})}_{+}/(\br^{+})^{d-n}$.
\end{lem}

\begin{proof} Let $\sigma=\langle v_1\ldots,v_n\rangle\in \Delta(n)$ be
such that $\{v_1,\ldots,v_n\}$ form a $\bz$ basis for $N$. Let
$\{u_1,\ldots,u_n\}$ be the dual basis. Let ${N}''\simeq \bz^{d-n}$;
$N'\simeq\bz^d$ and let $\{e_j':1\leq j\leq d\}$, $\{e_k'':1\leq k\leq
d-n\}$ denote the natural bases of $N'$ and $N''$
respectively. Further, let $g: N'\ra N$ map $e_j'$ to $v_j$ for every
$1\leq j\leq d$ and let $f:N''\hra N'$ be the map which takes
$e_{d-j+1}''$ to $e_j'-(\Sigma_{i=1}^ {n}\langle u_i,v_j\rangle\cdot
e_i')$ for every $n+1\leq j\leq d$.  From the results of \cite{Cox} we
know that there is an exact sequence of fans:
\begin{equation*}
0\lr (\Delta'', N'')\stackrel{f}\hra (\Delta',N')\stackrel{g}\lr
(\Delta,N)\lr 0,
\end{equation*}
where $\Delta'' =\{0\}$ and $\Delta'$ is the fan consisting of the
cones $\tau'=\langle e_{j_1}',\ldots,e_{j_k}'\rangle$ corresponding to
every $\tau=\langle v_{j_1},\ldots,v_{j_k}\rangle \in\Delta $. Observe
that, the real toric varieties corresponding to $\Delta''$ and $\Delta'$
are $ X(\Delta{''})\simeq (\br^*)^{d-n}$ and $X(\Delta')\simeq \br^d- Z$
respectively, where $Z$ is the zero locus in $\br^d$ of the monomials
$x_{\hat{\sigma}}=\prod_{\rho\not{\in}\sigma}x_{\rho}$ for every
$\sigma\in \Delta$. Moreover, it is easy to see that $\br^d- Z$ is also
isomorphic to the complement of the subspace arrangement
$\br^d-\ca_{\Delta}=\cc_{\Delta}$ defined above (cf. p.~130 of
\cite{Buck}). Hence, from the above exact sequence of fans, we see that
the smooth complete real toric variety $X$ is the base space of a
principal bundle with total space $(\br^d- Z)\simeq
(\br^d-\ca_{\Delta})\simeq \cc_{\Delta}$ and structure group
$(\br^*)^{d-n}$(cf. p.~59 of \cite{Oda} and p.~27 of \cite{Cox}).
Similarly, by restricting to the non-negative parts we see that $X_{+}$
is the base space of a principal fibre bundle with total space
$\br^d_{+}- Z_{+}$ and structure group $(\br^+)^{d-n}$. Thus we have the
following:
\begin{align*}
X &\simeq (\br^{d}- Z)/(\br^*)^{d-n}\simeq
\cc_{\Delta}/(\br^*)^{d-n},\\[.2pc]
X_{+} &\simeq \br_{+}^d- Z_{+}/(\br^+)^{d-n}\simeq
(\cc_{\Delta})_{+}/(\br^+)^{d-n}.
\end{align*}

$\left.\right.$\vspace{-2.5pc}

\hfill $\Box$
\end{proof}

\begin{rem}{\rm Note that the only property of a smooth and complete fan
which we use in the above proof is that $\{v_1,\ldots,v_n\}$ form a
$\bz$ basis of $N$. Thus Lemma~\ref{GQ} is true even for a smooth (not
necessarily complete) fan $\Delta$, for which the primitive vectors
along $\Delta(1)$ contains a $\bz$ basis for $N$.}
\end{rem}

\begin{lem}\label{fundsub} $\pi_1({\cc}_{\Delta})$ is isomorphic to the
commutator subgroup $[W,W]$ of the Coxeter group $W$ defined in {\rm \S2,}
which is generated as a normal subgroup of $W$ by $[s_{i_p}, s_{i_q}]$
for $1\leq i\leq r$ where{\rm ,} ${\cp}_i~=~\{\rho_{i_p},\rho_{i_q}\!\}~
\forall~ 1\leq i\leq r$.
\end{lem}

\begin{proof} From Lemma~\ref{GQ} we know that $X \simeq
\cc_{\Delta}/(\br^*)^{d-n}$. Moreover, since $(\br^*)^{d-n}\simeq
(\br^+)^{d-n}\times \bz_2^{d-n}$, \hbox{${X_{1}}=\cc_{\Delta}/(\br^+)^{d-n}$}
is a regular covering space over $X$ with deck transformation group
$\bz_2^{d-n}$. In fact, it is the same covering space of $X$ as in
Theorem~\ref{FG}(4). Also observe that $\cc_{\Delta}$ and ${X_{1}}$ are
of the same homotopy type since $\cc_{\Delta}$ is a fibre bundle over
${X_{1}}$ with contractible fibre $(\br^+)^{d-n}$. Therefore we have,
$\pi_1(\cc_{\Delta})\simeq[W,W]$.\hfill $\Box$
\end{proof}

In the following theorem we shall find the necessary and sufficient
conditions on $\Delta$ and hence on the arrangement $\ca_{\Delta}$,
under which $\cc_{\Delta}$ is aspherical.

\begin{theor}[\!]\label{aspher2}
$\cc_{\Delta}$ is aspherical if and only if $\ca_{\Delta}$ is a union of
precisely codimension $2$ subspaces.
\end{theor}

\begin{proof} Since $\cc_{\Delta}$ is of the homotopy type of a finite
regular covering space over $X$, it follows that $X$ is aspherical if
and only if $\cc_{\Delta}$ is aspherical. From Theorem~\ref{aspher1} the
necessary and sufficient condition for $X$ to be aspherical is that
${\Delta}$ is flag-like. Therefore, it suffices to show that ${\Delta}$
is flag-like if and only if $\ca_{\Delta}$ is a union of precisely
codimension two subspaces.

Now, by Definition~7.1, the condition for ${\Delta}$ to be flag-like
is equivalent to the condition that in $\Delta$ there are no primitive
collections consisting of more than two edges. Also by Definition~7.2,
$\ca_{\Delta}=\cup_{\cp}\ca_{\cp}$, where the union is over primitive
collections $\cp$ in $\Delta$ and where $\ca_{\cp}$ is a subspace in
$\br^d$ of codimension precisely equal to the number of edges in
$\cp$. Thus, ${\Delta}$ is flag-like if and only if
$\ca_{\Delta}=\cup_{\cp_{i}}\ca(\cp_{i})$ where the union runs over
the primitive collections $\{\cp_1,\ldots,\cp_r\}$ consisting of two
edges or equivalently, $\ca_{\Delta}$ is a union of codimension two
subspaces. Hence the theorem.\hfill $\Box$
\end{proof}

\begin{rem}\hskip -.5pc{\rm ({\it $K(\pi,1)$-arrangements}).\ \ The barycentric
subdivision of any simplicial complex is a flag complex. Hence, given a
smooth complete fan $\Delta$, we can obtain several smooth complete
flag-like fans whose cones are the cones over the simplices of the
repeated barycentric subdivisions of $\cs_{\Delta}$. We therefore get
several examples of $K(\pi,1)$ arrangements finding which seems to be of
interest in the topology of arrangements (cf. \cite{Orlik} and
\cite{khov}). However, note that even if we start with a flag-complex,
an arbitrary subdivision need not result in a flag-complex. For example,
let $\Delta$ be the fan consisting of the faces of $\sigma=\langle
e_1,e_2,e_3\rangle$ in $N=\bz e_1\oplus\bz e_2\oplus\bz e_3$. If we
refine $\Delta$ by adding the edge vector through $v=e_1+e_2+e_3$, then
the resulting fan $\Delta'$ is not flag-like since, $e_1,e_2,e_3,v $
pairwise span cones in $\Delta'$ but together do not span any cone.}
\end{rem}

\begin{rem}{\rm Indeed, both Lemma~\ref{fundsub} and Theorem~\ref{aspher2}
follow directly from the fact that $\cc_{\Delta}$ is a smooth
non-complete toric variety associated to the fan $\Delta'=\{\langle
e_{j_1},\ldots, e_{j_k}\rangle$ for every cone $\tau=\langle
v_{j_1},\ldots,v_{j_k}\rangle\in\Delta\}$ in $N'=\bz
e_1\oplus\cdots\oplus \bz e_d$ (cf. Lemma~\ref{GQ}) and applying
Theorems~\ref{FG} and \ref{aspher1}. However, since
$\cc_{\Delta}$ has been defined specifically as the complement of real
coordinate subspace arrangement related to a smooth complete fan
$\Delta$, we therefore describe both its fundamental group and
criterion for asphericity using $\Delta$.}
\end{rem}

\begin{rem}{\rm Since $\cc_{\Delta}$ is the toric variety associated to
the fan $\Delta'$, we can apply Theorems~\ref{pres} and \ref{ab}
respectively to give a presentation for $\pi_{1}(\cc_{\Delta})$ and give
conditions on $\Delta'$ for it to be abelian. In particular, it follows
from Theorems~\ref{aspher2} and \ref{ab} that $\cc_{\Delta}$ is
$K(\pi,1)$ with $\pi_1(\cc_{\Delta})$ abelian, if and only if it is the
complement of subspaces of codimension precisely $2$ which pairwise
intersect at $\{0\}$. Moreover, it also follows from Lemma~\ref{FO} that
$\pi_1(\cc_{\Delta})=[W,W]$ is always torsion free.}
\end{rem}

\section*{Acknowledgement}

I am grateful to Prof.~P~Sankaran for suggesting this problem, for his
invaluable guidance and constant encouragement in this work. I also
thank Prof.~V~Balaji for several helpful discussions. I thank the
referee for suggesting some corrections in the manuscript.


\begin{thebibliography}{99}
\bibitem{Audin} Audin~M, The topology of torus actions on
symplectic manifolds, {\it Prog. Math.} (Birkhauser, Basel, Boston
and Berlin) (1991) vol.~93

\bibitem{Brid} Bridson~M~R and Haefliger~A, Metric spaces of
non-positive curvature, A Series of Comprehensive Studies in
Mathematics (Springer) (1999) vol.~319

\bibitem{Brown} Brown~K~S, Buildings (Springer Verlag) (1999)

\bibitem{Buck} Buchstaber~V~M and Panov~T~E, Torus actions and their
applications in topology and combinatorics, {\it Univ. Lecture Series
AMS}, vol.~24 (2002)

\bibitem{Cohen} Cohen~D~E, Combinatorial group theory: A topological
approach, {\it London Math. Soc. Stud. Texts} {\bf 14} (1989)

\bibitem{Cox} Cox~D~A, The homogeneous coordinate ring of a toric
variety, {\it J. Algebraic Geometry} {\bf 4} (1995) 17--50

\bibitem{Davis} Davis~M~W, Groups generated by reflections and
aspherical manifolds not covered by the Eucledian space, {\it Ann.
Math.} {\bf 117(2)} (1983) 293--324

\bibitem{DJ} Davis~M~W and Januszkiewicz T, Convex polytopes Coxeter
orbifolds and torus actions, {\it Duke Math. J.} {\bf 62(2)} (1991)
417--451

\bibitem{DJS} Davis~M~W, Januszkiewicz~T and Scott~R, {\it Selecta Math.
(N.S.)} {\bf 4(4)} (1998) 491--547

\bibitem{Fult} Fulton~W, Introduction to toric varieties, {\it Ann.
Math. Studies} (Princeton University Press) (1993) no.~131

\bibitem{Jurk} Jurkiewicz~J, Torus embeddings, polyhedra, $k^*$-actions
and homology, {\it Dissertationes Math. (Rozprawny Mat.)} {\bf 236}
(1985) 1--69

\bibitem{khov} Khovanov~M, Real $K(\pi,1)$ arrangements from finite root
systems, {\it Math. Res. Lett.} {\bf 3} (1996) 261--274

\bibitem{LS} Lyndon~R and Schupp~P, Combinatorial group theory
(Springer) (1977)

\bibitem{Oda} Oda~T, Convex bodies and algebraic geometry, {\it
Ergebnisse der Mathematik} (Springer-Verlag) (1988) vol.~15

\bibitem{Orlik} Orlik~P, Arrangements in topology, discrete and
computational geometry (NJ: New Brunswick) (1989/1990) pp.~263--272,
DIMACS Ser. Discrete Math. Theoret. Comput. Sci., vol.~6, {\it Amer. Math.
Soc.} (Providence, RI) (1991)
\end{thebibliography}
\end{document}